\documentstyle[12pt]{article}
\input mssymb.tex
\ifx\msy\undefined\def\msy{\bf}\fi
\makeatletter                                            
\renewcommand{\@begintheorem}[2]{                        
\rm \trivlist \item [\hskip \labelsep {\bf #2\ \ #1.}]   
				}                        
\makeatother                                             

\makeatletter
\def\section{\@startsection {section}{1}{\z@}{-3.5ex plus -1ex minus
 -.2ex}{1.5ex plus .2ex}{\large\bf}}
\def\subsection{\@startsection{subsection}{2}{\z@}{-3.25ex plus -1ex
minus
 -.2ex}{-1em}{\normalsize\bf}}
\let\emppsubsection\subsection

\def\empsubsection[#1]#2{\emppsubsection[#1]{#2\unskip}}
\def\subsection{\secdef\empsubsection{\emppsubsection*}}

\makeatother

\newcommand{\newsubsubsection}%
{{\bf\refstepcounter{subsubsection}\thesubsubsection\ \ }}

\makeatletter
\let\c@equation\c@subsubsection
\makeatother

\newcommand{\CC}{{\msy C}}

\newcommand{\PP}{{\msy P}}
\newcommand{\QQ}{{\msy Q}}

\newcommand{\ZZ}{{\msy Z}}
\newcommand{\al}{\alpha}
\newcommand{\be}{\beta}

\newcommand{\G}{\Gamma}

\newcommand{\T}{\Theta}
\newcommand{\aje}{\approx_{AJ}}
\newcommand{\ale}{\approx_{alg}}
\newtheorem{Proposition}[subsection]{Proposition}
\newtheorem{Theorem}[subsection]{Theorem}
\newtheorem{Lemma}[subsection]{Lemma}

\newtheorem{Remark}[subsection]{Remark}
\newcommand{\qed}{{\unskip\nobreak\hfill\hbox{ $\Box$}\par}}

\begin{document}

\title{Note on curves in a Jacobian}
\author{Elisabetta Colombo and Bert van Geemen}
\date{\mbox{}}
\maketitle
\section{Introduction}
\subsection{}
For an abelian variety $A$ over $\CC$ and a cycle $\al\in CH_d(A)_\QQ$ we
define a subspace $Z_\al$ of $CH_d(A)_\QQ$ by:
$$
Z_\al:=\langle n_*\alpha:\;n\in\ZZ  \rangle \quad\subset CH_d(J(C))_{\QQ}.
$$
Results of Beauville imply that $Z_\al$ is finite dimensional (cf.\ \ref{short}
below). In case $A=J(C)$, the jacobian of a curve $C$, Ceresa has shown that
the cycle $C-C^-:=C-(-1)_*C\in Z_C$
is not algebraically equivalent to zero for generic $C$ of genus $g\geq 3$,
which implies that for such a curve $\dim_\QQ Z_C\geq 2$.

In this note we investigate the subspace $Z_{W_m}$ of $CH_m(J(C))_\QQ$, with
$W_m$ the image of the $m$-th symmetric power of $C$ in $J(C)$ (so $W_1=C$). To
simplify matters we will actually work modulo algebraic equivalence (rather
than linear equivalence, note that translates of a cycle are algebraically
equivalent).

\subsection{}
Let $Z_\al/\ale$ be the image of $Z_\al\subset CH_d(A)_\QQ$ in
$CH_d(A)_{\QQ}/\ale$.
A d-cycle $\al$ is Abel-Jacobi equivalent to zero, $\al\aje 0$, if $\al$ is
homologically equivalent to zero and its image in $J_d(A)$, the $d$-th
primitive Intermediate Jacobian of $A$, is zero. Recall that any curve of genus
$g$ is a cover of $\PP^1$ for some $d\leq\frac{g+3}{2}$.

\begin{Theorem}\label{thm}
\begin{enumerate}
\item For any abelian variety $A$ and any $\alpha\in CH_d(A)$ we have:
$$
\dim_\QQ \;(Z_\al/\aje)\leq 2.
$$
\item For any curve of genus $g$ and $1\leq n\leq g-1$ we have:
$$
 \dim_\QQ \; (Z_{W_{g-n}}/\ale)\leq n.
$$
\item For a curve $C$ which is a $d:1$-cover of $\PP^1$ we have:
$$
\dim_\QQ \; (Z_C/\ale)\leq d-1.
$$
\end{enumerate}
\end{Theorem}

(We prove \ref{thm}.1 in \ref{short}, \ref{thm}.2 in \ref{ct} and
\ref{thm}.3 in \ref{recurs}.)

\subsection{}
Recall that Ceresa showed that the image of $W_m-W_m^-$ in $J_m(J(C))$ is
non-zero for generic $C$ of genus $g\geq 3$ and $1\leq m\leq g-2$. Therefore
\ref{thm}.1 and \ref{thm}.2 for $n=2$ are sharp.
In case $C$ is hyperelliptic, so $C$ is a 2:1 cover of $\PP^1$, the cycles
$W_m$ and $W_m^-$ are however algebraically equivalent. Therefore \ref{thm}.3
is sharp for hyperelliptic curves ($d=2$) and generic trigonal curves ($d=3$).
In case $C$ is not hyperelliptic nor trigonal (like the generic curve of genus
$g\geq 5$), it would be interesting to know if \ref{thm}.3 is actually sharp.

Note that \ref{thm}.1 implies that one cannot use the Intermediate Jacobian
anymore to derive new algebraic relations among the cycles in $Z_C/\ale$.
Recently M.Nori \cite{N} constructed cycles on complete intersections in
$\PP^N$ which are Abel-Jacobi equivalent to zero but not algebraically
equivalent to zero. There is thus the possibility that similar cycles can be
found on the Jacobian of a curve of genus $g\geq 5$.
A cycle in $Z_C$, for certain modular curves $C$, was investigated by B.H.Gross
and C.Schoen, \cite{GS}, esp.\ section 5, see also \ref{sgz}.

\subsection{} The inequalities \ref{thm}.1 and \ref{thm}.2 are consequences of
work of Beauville.
The main part of the paper deals with the proof of \ref{thm}.3. Recall that on
a smooth surface $S$ homological and algebraical equivalence for curves
coincide. Thus if we have map $\Phi:S\rightarrow J(C)$ and relation
$a_1[C_1]+\ldots +a_n[C_n]=0$ in $H^2(S,\QQ)$, we get $a_1\Phi_*C_1+\ldots
+a_n\Phi_*C_n\ale 0$ in $J(C)$. We use this remark to obtain our result, the
main difficulty is of course to find suitable surfaces, curves in them and maps
to $J(C)$ and to determine Neron-Severi groups ($=Im(CH_1(S)\rightarrow
H^2(S,\QQ))$) of the surfaces involved.

\subsection{} We are indebted to S.J.Edixhoven and C.Schoen for several helpful
discussions.

\section{General results}

\subsection{}
The effect of $n_*$ and $n^*$ on the Chow groups has been investigated by
Beauville (\cite{B1}, \cite{B2}), Dehninger and Murre \cite{DM} and K\"unnemann
\cite{K} using the Fourier transform on abelian varieties. Below we summarize
some of their results and derive the finite dimensionality of $Z_\alpha$ as
well as Thm \ref{thm}.1 and \ref{thm}.2.

\subsection{}
 Let $A$ be a $g$-dimensional abelian variety, we will view $B:=A\times A$ as
 an abelian scheme over $A$ using the projection on the first factor:
{\renewcommand{\arraystretch}{1.5}
$$
\begin{array}{c} B\\ \downarrow \\ A \end{array} \quad =\quad
\begin{array}{l} A\times A \\ \downarrow \pi_1 \\ A \end{array}
$$
}
For each integer $n$, we have an $s_n\in B(A)$ and a cycle $\Gamma_n\in
CH^g(B)$:
$$
s_n:A\longrightarrow B,\quad a\longmapsto (a,na),\qquad
\Gamma_n:=s_{n*} A\in CH^g(A\times A),
$$
and $\Gamma_n$ is, as the notation suggests, the graph of multiplication by $n$
on $A$. The cycle $\Gamma_n$ defines for each $i$ a map on $CH^i(A)$ which is
just $n_*$:
$$
n_*=\Gamma_n:CH^i(A)\longrightarrow CH^i(A),\quad
\al\longmapsto n_*\al=\pi_{2*}(\pi_1^*\al\cdot \Gamma_n).
$$

Next we introduce a relative version of the Pontryagin product $*$ on the Chow
groups of the $A$-scheme $B$ (cf. \cite{K}, (1.2)). Let $m_B:B\times_A
B\rightarrow B$ be the multiplication map, then:
$$
\al*\be:=m_{B*}(\al\times_A \be),\qquad \al,\be\in CH^*(B).
$$
Since the relative dimension of $\Gamma_n$ over $A$ is $0$, the cycle
$\Gamma_n*\Gamma_m$ lies in  $CH^g(B)$ and one has (\cite{K}, (1.3.4)):
$$
\Gamma_n * \Gamma_m =\Gamma_{n+m}\quad \in CH^g(A\times A)_\QQ,\qquad{\rm
so}\quad
\Gamma_n=\Gamma_1^{*n},
$$
where we write $\al^{*n}$ for the $n$-fold Pontryagin product of a cycle $\al$
with $n>0$ and we put $\al^{*0}:=\Gamma_0$. In \cite{K} 1.4.1 a generalization
of a theorem of Bloch is proved, which implies:
\begin{equation}\label{BG}
(\Gamma_1-\Gamma_0)^{*(2g+1)}=0\qquad (\in CH^g(B)_\QQ).
\end{equation}
Using the ring structure on $CH^g(B)_\QQ$ with product $*$ one can thus define
the following cycles $\pi_i$, $0\leq i\leq 2g$ in $CH^g(B)_\QQ$:
$$
\pi_i:=\mbox{$\frac{1}{(2g-i)!}$}(\log \G_1)^{*(2g-i)},
$$
with:
$$
\log \G_1:=(\G_1-\G_0)-\mbox{$\frac{1}{2}$}(\G_1-\G_0)^{*2}+\ldots +
\mbox{$\frac{1}{2g}$}(\G_1-\G_0)^{*(2g)}.
$$
Let $\Delta=\G_1\in CH^g(A\times A)_\QQ$ the class of the diagonal, then
(\cite{K}, proof of 3.1.1):
\begin{equation}\label{del}
\Delta=\pi_0+\pi_1+\ldots + \pi_{2g},\qquad
\pi_i\pi_j=\pi_j\pi_i=\left\{\begin{array}{ll} \pi_i&{\rm if}\;\;i=j\\
						0&{\rm if}\;\;i\neq j,
			    \end{array}\right.
\end{equation}
here $\al\be$, for cycles $\al, \be\in CH^g(A\times A)_\QQ$, is their product
as correspondences: $\al\be:=p_{13*}(p_{12}^*\al\cdot p_{23}^*\be)$ with the
$p_{ij}:A^3\rightarrow A^2$ the projection to the $i,j$ factor. Moreover:
\begin{equation}\label{gp}
\G_n\pi_{2g-i}=\pi_{2g-i}\G_n=n^i\pi_{2g-i},
\end{equation}
(one has ${}^t\G_n\pi_i=\pi_i{}^t\G_n=n^i\pi_i$ and ${}^t\pi_i=\pi_{2g-i}$
(\cite{K}, 3.1.1), now take transposes).

\begin{Remark} We sketch how these results can be obtained from \ref{BG}.
Let $M\subset CH^g(B)_\QQ$ be the subspace spanned by the $\G_n$, $n\in\ZZ$.
Then \ref{BG} implies that $\dim_\QQ\, M\leq 2g+1$ (use
$\G_i*(\Gamma_1-\Gamma_0)^{*(2g+1)}=0$ for all $i\in \ZZ$).
Using the K\"{u}nneth formula, Poincar\'e duality,
one finds that the cohomology class of $\G_n$ in
$$
H^{2g}(A\times A,\QQ)=\oplus H^{2g-i}(A,\QQ)\otimes H^i(A,\QQ)=\oplus
{\rm Hom}(H^i(A,\QQ),H^i(A,\QQ))
$$
is given by (note $\G_n$ induces $n_*$):
$$
[\G_n]=(n^{2g},n^{2g-1},\ldots ,n,1)\in \oplus_{i=0}^{2g}
{\rm Hom}(H^i(A,\QQ),H^i(A,\QQ))
$$
(where $n^{2g-i}$ in the $i$-th component means multiplication by $n^{2g-i}$ on
$H^i(A,\QQ)$). Therefore $\dim_\QQ\,M/\sim_{hom} =2g+1$ and $\dim_\QQ\, M=
2g+1$. Thus we have the (surprising) result that homological and linear
equivalence coincide in $M$. We can now define $\pi_i\in M$ by
$$
[\pi_i]:=(0,\ldots ,0,1,0,\ldots ,0)\in \;\oplus_i
{\rm Hom}(H^i(A,\QQ),H^i(A,\QQ))
$$
(1 in $i$-th spot) then \ref{del} and \ref{gp} follow. To express $\pi_i$ as
combination of
$\G_0,\ldots ,\G_{2g}$, note that the ($\QQ$-linear) ring homomorphism:
$$
\QQ[X]\longrightarrow M\subset CH^g(B)_\QQ,\qquad X^i\longmapsto \G_i=\G_1^{*i}
$$
(with $*$ product on $CH^g(B)_\QQ$) gives an isomorphism
$\QQ[X]/(X-1)^{2g+1}\cong M$. Since $\G_n\G_m=\G_{nm}$ (product as
correspondences), $\pi_{2g-i}$ corresponds to a polynomial $f_{2g-i}$ with:
$$
f_{2g-i}(X^n)=n^if_{2g-i}(X)\qquad{\rm so}\quad
f_{2g-i}(X):=c_i(\log X)^i\;{\rm mod}\;(X-1)^{2g+1}
$$
and $c_i\in\QQ$ can be determined with a little more work.
\end{Remark}

\subsection{} Since $\Delta=\G_1:CH_d(A)\rightarrow CH_d(A)$ is the identity,
we get:
$$
CH_d(A)_\QQ=\oplus CH_d(A)_{(i)},\qquad{\rm with}\quad
 CH_d(A)_{(i)}:=\pi_{2g-i}CH_d(A)_\QQ.
$$
and each $CH_d(A)_{(i)}$ is an eigenspace for the multiplication operators:
$$
n_*\al=n^i\al\qquad\forall n\in \ZZ,\;\;\forall\al \in CH_d(A)_{(i)},
$$
a result which was first obtained by Beauville \cite{B2}. Moreover, he proves:
\begin{equation}\label{bb}
CH_d(A)_{(i)}\neq 0\qquad \Longrightarrow \quad d\leq i\leq d+g
\end{equation}
(and gives sharper bounds for some $d$ in prop.3 of \cite{B2}, it has been
conjectured that $CH_d(A)_{(i)}\neq 0\;\Leftrightarrow\; 2d\leq i\leq d+g$).

\begin{Proposition}\label{short}
 Let A be an abelian variety and let $\al\in CH_d(A)_\QQ$. Then:
$$
\dim_\QQ Z_\al \leq g+1\qquad {\rm and}\quad \dim_\QQ (Z_\al/\aje) \leq 2.
$$
Moreover, we have $\pi_i Z_\al \subset Z_\al$ for all $i$.
\end{Proposition}
{\bf Proof.}$\quad$We write $\al$ as a sum of weight vectors
($\al_{(i)}:=\pi_{2g-i}\al$):
$$
\al=\al_{(d)}+\al_{(d+1)}+\ldots +\al_{(g+d)},\qquad {\rm so}\quad
n_*\al_{(i)}=n^i\al.
$$
Taking $g$ distinct, non-zero, integers $n_j$, the determinant of the matrix
expressing the $n_{i*}\al\in Z_\al$ in terms of the $\al_i$ is a Vandermonde
determinant. Thus each $\al_{(i)}=\pi_{2g-i}\al\in Z_\al$ and $Z_\al$ is
spanned by the $\al_{(i)}$.

Since $n_*$ acts as $n^{2d}$ on $H^{2g-2d}(A,\QQ)$ and as $n^{2d+1}$ on
$J_d(A)$, the space $Z_\al/\aje$ is spanned by $\al_{(2d)}$ and
$\al_{(2d+1)}$.

\qed

\subsection{} In this section, we fix an abel-jacobi map $C\hookrightarrow
J(C)$. Note that $W_d=\mbox{$\frac{1}{d!}$}C^{*d}$ (Pontryagin product on
$J(C)$). Let $\Theta\in CH^1(J(C))$ be a symmetric theta divisor, so $\T$ is a
translate of $W_{g-1}$. Let $\Theta^d\in CH_{g-d}(J(C))$ be the $d$-fold
intersection of $\T$.

\begin{Proposition}\label{ct}
 For any $d,\;1\leq d\leq g-1$, we have $\Theta^d\in Z_{W_{g-d}}$, more
 precisely:
$$
\pi_{2g-i}W_{g-d}\neq 0\;\Longrightarrow \; 2(g-d)\leq i\leq 2g-d,
\quad{\rm and}\quad
\pi_{2d} W_{g-d}= \mbox{$\frac{1}{d!}$}\Theta^d.
$$
Moreover we have:
$$
\dim_\QQ \;(Z_{W_{g-d}}/\ale) \leq d.
$$
\end{Proposition}
{\bf Proof.}$\quad$We first prove the case $d=g-1$. Since the map
$CH_1(C)_\QQ\rightarrow H^{2g-2}(J(C),\QQ)$ factors over
$\pi_{2g-2}CH_1(J(C))_\QQ$ (cf.\ the proof of \ref{short} or even better,
\cite{B2}, p.650) we know that $\pi_{2g-2}C\neq 0$. Then its Fourier transform
$F_{CH}(\pi_{2g-2}C)\in \pi_2CH_{g-1}(J(C))_\QQ\;(\cong NS(J(C))_\QQ)$ (cf.
\cite{B2}, prop.1) is non-zero.
For a generic Jacobian $NS_\QQ$ is one dimensional and thus
$\pi_2CH_{g-1}(J(C))_\QQ=\QQ\T$ (\cite{B1}, prop.5 and \cite{B2}, prop.1).
By specializing, we have for all curves that $F_{CH}(\pi_{2g-2}C)\in \QQ\T$.

 From \cite{B1}, prop.5 we have
$\Theta=-F_{CH}(\mbox{$\frac{1}{(g-1)!}$}\Theta^{g-1})$,
so for a nonzero constant $c$:
$$
F_{CH}(\pi_{2g-2}C)=cF_{CH}(\mbox{$\frac{1}{(g-1)!}$}\T^{g-1})$$
and, by using $F_{CH}^2=(-1)^g\G_{-1}$ and comparing cohomology classes, we
get:
$$
    \pi_{2g-2}C=\mbox{$\frac{1}{(g-1)!}$}\T^{g-1}.
$$

Next we recall that $\pi_{2g-1}CH_1(A)=0$ for any abelian variety $A$
(\cite{B2}, prop.3), thus:
$$
C=C_{(2)}+\ldots + C_{g+1} \qquad {\rm with~} C_{(i)}=\pi_{2g-i}C,
$$
so $n_*C_{(i)}=n^iC_{(i)}$ and $C_{(2)}=\mbox{$\frac{1}{(g-1)!}$}\T^{g-1}$.
Therefore
\begin{equation} \label{W}
W_{g-d}=\mbox{$\frac{1}{(g-d)!}$}C^{*(g-d)}=
\mbox{$\frac{1}{(g-d)!}$}C_{(2)}^{*(g-d)}+Y
\end{equation}
with $Y$ a sum of cycles $C_{(i_1)}*\ldots *C_{(i_{g-d})}$ with all $i_j\geq 2$
and at least one $>2$. Therefore, using
$n_*(U*V)=(n_*U)*(n_*V)$, we get:
$\pi_{2d}W_{g-d}=\mbox{$\frac{1}{(g-d)!}$}C_{(2)}^{*(g-d)}$ and
$\pi_jW_{g-d}=0$ for $j>2d$.

By \cite{B1}, corr.2 of prop.5, $\Theta^d$ lies in the subspace spanned by
$(\Theta^{g-1})^{*(g-d)}$, so $C_{(2)}^{*(g-d)}=c\Theta^d$  for a non-zero
$c\in\QQ$ and taking cohomology classes one finds $c=1$.

Finally, using \ref{bb} and \ref{W} we can write
$$
W_{g-d}=\alpha_{(2g-2d)}+\ldots +\alpha_{(2g-d)},\qquad{\rm with}\quad
\al_{(i)}:=\pi_{2g-i}W_{g-d},
$$
and since $\al_{(2g-d)}\ale 0$ by \cite{B2}, prop.4a, we get $\dim_\QQ
(Z_{W_{g-d}}/\ale) \leq d$.
\qed

\subsection{}\label{sgz}
In the paper \cite{GS} the following cycle (modulo $\ale$) is considered:
$$
Z:=3_*C-3\cdot 2_*C +3C = (\G_1-\G_0)^{*3}C.
$$
This cycle is related to Ceresa's cycle in the following way:

\begin{Proposition}
We have:
$$
\pi_{2g-2}(Z)=\pi_{2g-2}(C-C^-)=0\qquad{\rm and}\quad
\pi_{2g-3}(Z)=3\pi_{2g-3}(C-C^-).
$$
In particular, $Z$ is abel-jacobi equivalent to $3(C-C^-)$.
\end{Proposition}
{\bf Proof.}$\quad$As we saw before (proof of \ref{ct}) we can write:
$$
C=C_{(2)}+C_{(3)}+\ldots \qquad {\rm with}\;C_{(i)}:=\pi_{2g-i}C.
$$
Since $n_*C_{(i)}=n^iC_{(i)}$ we have:
$$
C-C^-=2( C_{(3)}+0+C_{(5)}+\ldots ),\qquad Z=6C_{(3)}+36C_{(4)}+\ldots.
$$
Therefore $\pi_{2g-2}(C-C^-)=3\pi_{2g-2}(Z)=0$, so both cycles are
homologically equivalent to zero, and
$3\pi_{2g-3}(C-C^-)=\pi_{2g-3}(Z)=6C_{(3)}$. Since the abel-jacobi map factors
over $\pi_{2g-3}CH_1(J(C))$ we have that
$Z\aje 3(C-C^-)$.
\qed

\section{Proof of 1.3.3}
\subsection{}
Let $C$ be a generic $d$:1 cover of $\PP^1$, and denote by $g^1_d$ the
corresponding linear series.
For an integer $n$, $0<n<d$ we define a curve $G_n$ in the $n$-fold
$C^{(n)}$=Sym$^n(C)$:
$$
G_n=G_n(g^1_d)=\{x_1+x_2+\ldots+x_n\in C^{(n)}:\quad
x_1+x_2+\ldots+x_n<g^1_d\},
$$
see \cite{ACGH}, p.342 for the definition of the scheme-structure on $G_n$. We
define a surface by:
$$
S_n=S_n(g^1_d)=G_n\times C.
$$

\begin{Proposition}\label{NS}
Let $C$ be a generic $d$:1-covering of $\PP^1$ of genus $g\geq 1$.

Then $S_n$ is a smooth, irreducible surface and
$$
\dim_\QQ NS(S_n)_\QQ=3.
$$
\end{Proposition}
{\bf Proof.}$\quad$Since $C$ is generic, the curve $G_n$ is irreducible
(consider the monodromy representation) and for a simple covering the
smoothness of $G_n$ follows from a local computation. Therefore also $S_n$ is
smooth and irreducible.

Note that for a product of 2 curves $S_n=G_n\times C$ we have that
$$
NS(S_n)_\QQ=\QQ C\oplus \QQ G_n\oplus {\rm Hom}(J(C),J(G_n))_\QQ.
$$
Using Hodge structures, the last part are just the Hodge-cycles in
$H^1(C,\QQ)\otimes H^1(G_n,\QQ)$.
Since $J(C)$ and $Pic^0(C^{(n)})$ are isogeneous (their $H^1(\QQ)$'s are
isomorphic Hodge structures), composing such an isogeny with the pull-back
$Pic^0(C^{(n)})$ $\rightarrow$ $J(G_n)$, we have a non-trivial map
$J(C)\rightarrow J(G_n)$. To show that $\dim_\QQ {\rm Hom}(J(C),J(G_n))_\QQ\leq
1$ we use a degeneration argument. First we recall that the Jacobian of $C$ is
simple and has in fact $End (J(C))=\ZZ$.

In case $g=2$ this is clear since any genus 2 curve can be obtained by
deformation from a given $d$:1 cover and $C$ is generic. In case $g=3$ one
reasons similarly, taking the necessary care for the hyperelliptic curves.
In case $g(C)=g>3$, assume first that there is an elliptic curve in $J(C)$.
Specializing $C$ to a reducible curve with two components $C'$ and $C''$, both
of genus $\geq 2$, and themselves generic $d$:1-covers we obtain a
contradiction by induction. Assume now that there is no elliptic curve in
$J(C)$, then we specialize $C$ to $E\times C'$, with $E$ an elliptic curve and
$C'$ a generic $d$:1-covering of genus $g-1$. By induction again, $J(C')$ is
simple, which contradicts the existence of abelian subvarieties in $J(C)$. We
conclude that $J(C)$ is simple. Thus $End_\QQ (J(C))$ is a division ring and
specializing again we find it must be $\QQ$.

We may assume that the $g^1_d$ exhibits $C$ as a simple cover of $\PP^1$, then
$G_n$ is a cover of $\PP^1$ (of degree $({}^d_n)$) with only twofold
ramification points.
Letting two branch points coincide, we obtain a curve $\overline{C}$ with a
node, the normalization $C'$ of $ \overline{C}$ has genus $g-1$ and is again
exhibited as a generic $d:1$ cover of $\PP^1$.
The curve $G_n$ acquires $({}^{d-2}_{n-1})$ nodes (since twice that number of
branch points coincide pairwise) and the normalization of that curve,
$\overline{G}_n$, is $G'_n$, the curve obtained from the $g^1_d$ on $C'$.
The (generalized) Jacobians of $\overline{C},\,\overline{G}_n$ are extensions
of the abelian varieties (of $\dim \geq 1$) $J(C'),\,J(G'_n)$ by multiplicative
groups.

Since the number of simple factors of $J(G_n')$ which are isogeneous to $J(C')$
is greater then or equal to the number of simple factors of $J(G_n)$ which are
isogeneous to $J(C)$, and $J(C),\, J(C')$ have $End_\QQ=\QQ$ it follows:
$$
\dim_\QQ {\rm Hom}(J(C),J(G_n))_\QQ\; \leq \dim_\QQ \;
{\rm Hom}(J(C'),J(G'_n))_\QQ.
$$
Therefore it suffices to show that for a generic elliptic curve $C=E$ and a
generic $g^1_d$ on $E$ we have $\dim_\QQ {\rm Hom}(E,J(E_n))_\QQ\leq 1$ (with
$E_n=G_n(g^1_d)$ and $C=E$).

We argue again by induction, using degeneration. Note that for any $C$,
$G_1(g^1_d)\cong G_{d-1}(g^1_d)\cong C$. Since a generic elliptic curve has
${\rm Hom}(E,E)=\ZZ$ the statement is true for $d\leq 3$ (and any $0<n<d$) and
since $E_1\cong E$ it is also true for $n=1$ (and any $d\geq 2$).

We fix $E$, a generic elliptic curve, but let the $g^1_d$ acquire a base point
$y\in E$. Then $E_n(g^1_d)$, $n\geq 2$, becomes a reducible curve $\bar{E_n}$,
having two components $E_n'=E_n(g^1_{d-1})$ and $E'_{n-1}=E_{n-1}(g^1_{d-1})$
which meet transversally in $N=({}^{d-2}_{n-1})$ points. Indeed, let the
divisor of the $g^1_{d-1}$ containing $y$ be $y+y_1+\ldots +y_{d-2}$, then the
points $y_{i_1}+\ldots +y_{i_{n-1}}\in E'_{n-1}$ and
$y+y_{i_1}+\ldots +y_{i_{n-1}}\in E'_{n}$ are identified.

Thus $J(\overline{E}_n)$ is an extension of $J(E'_n)\times J(E'_{n-1})$ by the
multiplicative group $(\CC^*)^{N-1}$, in fact there is an exact sequence:
$$
1\longrightarrow \CC^*\stackrel{\Delta}{\longrightarrow} (\CC^*)^N
\longrightarrow J(\overline{E}_n) \stackrel{\pi}{\longrightarrow} J(E'_n)\times
J(E'_{n-1}) \longrightarrow 0,
$$
where $\Delta$ is the diagonal embedding. We will write $\pi=(\pi_0,\pi_1)$.

We will prove that
$$
 {\rm Hom}(E,J(\overline{E}_n))_\QQ \longrightarrow {\rm Hom}(E,J(E'_n))_\QQ,
\qquad\tilde\phi\mapsto \pi_0\circ \tilde\phi
$$
is injective. By induction we may assume that $\dim{\rm Hom}(E,J(E'_n))_\QQ
\leq 1$, thus also $\dim_\QQ {\rm Hom}(E,J(\overline{E}_n))_\QQ\leq 1$ and
since $\dim_\QQ {\rm Hom}(E,J({E}_n))_\QQ $  $\leq \dim_\QQ {\rm Hom}
(E,J(\overline{E}_n))_\QQ$ the assertion on the rank of the Neron-Severi group
follows.

Assume that $\tilde\phi\neq 0$, but $\pi_0\circ\tilde\phi=0$. Then
$\tilde\phi(E)\subset \bar{J}:=\pi_1^{-1}(J(E'_{n-1}))$. Pulling back this
 $(\CC^*)^{N-1}$-bundle $\bar{J}$ over $J(E_{n-1}')$ to $E$ along
 $\pi_1\circ\tilde\phi$, we obtain a $(\CC^*)^{N-1}$-bundle $\tilde{E}$ over
 $E$.
The map $\tilde\phi:E\rightarrow \bar{J}$ gives a section of
$\tilde{E}\rightarrow E$, thus $\tilde{E}$ is a trivial
$(\CC^*)^{N-1}$-bundle.
We show that this gives the desired contradiction.

{\renewcommand{\arraystretch}{1.5}
\vspace{\baselineskip}
$$
\begin{array}{cccccccc}
1&\rightarrow & (\CC^*)^{N-1} & \rightarrow &\bar{J}&\rightarrow &
J(E'_{n-1})&\rightarrow 0\\
\uparrow&     &\uparrow      &              &\uparrow &           &
\uparrow \pi_1\tilde\phi&  \\
1&\rightarrow &(\CC^*)^{N-1}&\rightarrow &\tilde{E}&\rightarrow &
 E&\rightarrow 0
\end{array}
$$
}

For any distinct $P,Q\in E'_n$ which are in $E'_n\cap E'_{n-1}$, the
$(\CC^*)^{N-1}$-bundle $\bar{J}$ over $J(E_{n-1}')$ has a quotient
$\CC^*$-bundle $\bar{J}_{PQ}$ whose extension class is $P-Q\in
Pic^0(E'_n)=Ext^1(J(E'_n),\CC^*)$. By induction, there is a `unique' map in
${\rm Hom}(J(E'_{n-1}),E)$ which must thus be induced by $x_1+\ldots
+x_{n-1}\mapsto x_1+\ldots +x_{n-1}-(n-1)p,\; E'_{n-1}\rightarrow E$ for some
$p\in E$. Taking $P=y_1+\ldots +y_{n-1}$ and $Q=y_2+\ldots +y_{n-1}+y_n$, the
pull-back of $P-Q$ to $E=Pic^0(E)$ is $y_1-y_n$. Choosing the degeneration
suitably we may assume that $y_1-y_n$ is not a torsion point on $E$ and thus
$\bar{J}_{PQ}$ has a nontrivial pull-back to $E$, contradicting the fact that
$\tilde{E}$, the pull-back of $\bar{J}$ to $E$, is trivial.
\qed

\subsection{}
We define a curve $\Delta_n$ in the surface $S_n=G_n\times C$:
$$
\Delta_n:=\{(x_1+x_2+\ldots+x_n,p)\in S_n:\;p\in\{x_1,\ldots ,x_n\} \}
$$
and for $n+1\leq d$, another curve $H_n$ in $S_n$:
$$
H_n=\{(x_1+x_2\ldots+x_n,p)\in S_n:\;\; x_1+x_2+\ldots+x_n+p<g^1_d\}.
$$
Finally we define the map:
\begin{equation}
\Phi^n_{l,k}:S_n\longrightarrow J(C),
\end{equation}
$$
(x_1+x_2\ldots+x_n,p)\longmapsto l(x_1+x_2+\ldots+x_n)+kp-D_{nl+k},
$$
where $D_{nl+k}$ is some divisor of degree $nl+k$ on $C$. Finally we denote by
$$
 u_n:C^{(n)}\longrightarrow J(C),\qquad D\longmapsto D_n,
$$ with $D_n$ some divisor of degree $n$, the Abel-Jacobi map on $C^{(n)}$. We
simply write $C$ for $u_{1*}C$.

\begin{Proposition} Let $C$ be a generic $d$:1 cover of $\PP^1$ of genus $g\geq
1$. Then:
\begin{enumerate}
\item
$$
NS(S_n)_\QQ=\langle C,G_n,\Delta_n\rangle.
$$
\item In $NS(S_n)$ we have:
\begin{equation}\label{eq2}
H_n=({}^d_{n})C+dG_{n}-\Delta_n.
\end{equation}

 \item The image of $G_n$ in $J(C)$ is a combination of $C$, $2_*C,\ldots
 ,n_*C$:
\begin{equation}\label{eq3}
u_{n*}G_n\ale ({}^d_{n-1})C-\mbox{$\frac{1}{2}$}({}^d_{n-2})2_*C+\ldots+
\mbox{$\frac{(-1)^{n-1}}{n}$}({}^d_0)n_*C.
\end{equation}

\end{enumerate}

\end{Proposition}
{\bf Proof.}$\quad$
The first part follows from the previous proposition and from the fact that in
the proof of \ref{eq2} we will see that $H_n$ can be uniquely expressed as a
combination of $C$, $G_n$ and $\Delta_n$.

 The proof of 2 and 3 is by induction. In the case $n=1$, so $G_1\cong C$,
 \ref{eq3} is trivial and for \ref{eq2} we have $S_1=C\times C$ and we must
 show $H_1=d(C+G_1)-\Delta$. Note the following intersection numbers:
$$
H_1\cdot C=H_1\cdot G_1=d-1,\quad H_1\cdot \Delta=2(g+d-1)
$$
(use that the $g^1_d$ exhibits $C$ as $d$:1 cover of $\PP^1$ with $2(g+d-1)$
simple ramification points), which imply the result.

Suppose that \ref{eq2} and \ref{eq3} are true for all $k<n$.
Note that $u_{n}(G_n)=\Phi^{n-1}_{1,1}(H_{n-1})$ and that map
$\Phi^{n-1}_{1,1}$ restricted to $H_{n-1}$ is $n$:1 so, using \ref{eq2} for
$n-1$, we have in $J(C)$:
\begin{equation}\label{eq4}
G_n\ale  \mbox{$\frac{1}{n}$}\left(
({}^d_{n-1})C+dG_{n-1}-(\Phi^{n-1}_{1,1})_*\Delta_{n-1}\right).
\end{equation}
where we write $G_k$ for $u_{k*}G_k$.
Next we observe that
\begin{equation}\label{fies}
(\Phi^{k}_{1,l})_*\Delta_k=(\Phi^{k-1}_{1,l+1})_*H_{k-1}
\qquad{\rm and}\qquad
(\Phi^{k-1}_{1,l+1})_*C=(l+1)_*C.
\end{equation}
Using \ref{eq2} in the cases $n-2,\,n-3,\ldots, 1$ we get:
\begin{eqnarray}\nonumber
(\Phi^{n-1}_{1,1})_*\Delta_{n-1}& \ale &(\Phi^{n-2}_{1,2})_*H_{n-2}\\ \nonumber
&\ale & ({}^d_{n-2})2_*C+dG_{n-2}-(\Phi^{n-3}_{1,3})_*H_{n-3} \\ \label{eq5}
&\ale & ({}^d_{n-2})2_*C-({}^d_{n-3})3_*C+
d(G_{n-2}-G_{n-3})+(\Phi^{n-4}_{1,4})_*H_{n-4}\\ \nonumber
&\ale & ({}^d_{n-2})2_*C-({}^d_{n-3})3_*C+\ldots+
	(-1)^{n-2}({}^d_{0})n_*C+\\ \nonumber
& &+d(G_{n-2}-G_{n-3}\ldots+(-1)^{n-2}G_2+(-1)^{n-1}G_1)
\end{eqnarray}
(note that $(\Phi^{1}_{1,n-1})_*H_{1}\ale({}^d_1)(n-1)_*C+dG_1-({}^d_0)n_*C$).
Substitute the expression for $(\Phi^{n-1}_{1,1})_*\Delta_{n-1}$ from \ref{eq5}
into \ref{eq4}:

\begin{eqnarray} \nonumber
\;\;G_n &\hbox{\rlap{$\ale$}}& \\ \nonumber
& & \!\mbox{$\frac{1}{n}$}\left(({}^d_{n-1})C-({}^d_{n-2})2_*C\ldots +
(-1)^{n-2}({}^d_{1})(n-1)_*C+(-1)^{n-1}({}^d_{0})n_*C \right. \\ \label{Gn}
 & &+\left. d(G_{n-1}-G_{n-2}+\ldots +(-1)^{n-3}G_2+(-1)^{n-2}G_1\right).
\end{eqnarray}

Using the formula \ref{eq3} for $G_k$, $k=1,2,\ldots ,n-1$ and the relation
$\sum_{k=0}^l(-1)^k({}^d_{l-k})=({}^{d-1}_l)$, (cf. \ref{bin1})
we get:
\begin{eqnarray}\label{Gs}
& & G_{n-1}-G_{n-2}+\ldots+(-1)^{n-2}G_1 \\ \nonumber
&\ale &({}^{d-1}_{n-3})C-\mbox{$\frac{1}{2}$}({}^{d-1}_{n-4})2_*C+\ldots+
(-1)^{n-3}\mbox{$\frac{1}{n-1}$}({}^{d-1}_0)(n-1)_*C.
\end{eqnarray}
Substituting this in formula \ref{Gn} and using the identity:
$$
\mbox{$\frac{1}{n}$}[({}^d_{n-k})+\mbox{$\frac{d}{k}$}({}^{d-1}_{n-k+1})]=
\mbox{$\frac{1}{k}$}({}^d_{n-k})
$$
we obtain \ref{eq3}.

To obtain \ref{eq2}, note that by (1), there are $a,b,c\in \QQ$ such that:
\begin{equation}\label{eq6}
H_{n}=aC+bG_n+c\Delta_n.
\end{equation}
To find them we compute the homology classes of the curves in \ref{eq6} in
$H_2(J(C),\QQ)$ and the intersection numbers of $H_n$ with $C$ and $G_n$.

For the generic curve $C$, the homology class $[B]$ of a curve $B$ in $J(C)$ is
a multiple of the class $[C]$ of $C$ and this multiple is
$\mbox{$\frac{1}{g}$}\T\cdot B$. We apply this to $B=G_n$. Using the map
$u:C^{(n)}\rightarrow J(C)$, we have $u_*(G_n)\cdot \T=u_*(G_n\cdot \theta)$
with $\theta$ the pull-back of $\T$ to $C^{(n)}$.

The homology class of $G_n$ in $H_2(C^{(n)},\QQ)$ is:
\begin{equation}
G_n=\sum_{k=0}^{n-1}
({}^{d-g-1}_k)\frac{x^k\theta^{n-k-1}}{(n-1-k)!},\qquad{\rm so}\quad
G_n\cdot\theta=\sum_{k=0}^{n-1} ({}^{d-g-1}_k)\frac{x^k\theta^{n-k}}{(n-1-k)!},
\end{equation}
cf.\ formula (3.2) on pag.342 of \cite{ACGH} (here $x$ is the class of the
divisor $C^{(n-1)}$ in $C^{(n)}$). From pag.343 of \cite{ACGH} one has:
$ u_*(x^{n-i}\theta^i)=({}^g_i)i![\T^g]/g! $
and since $[\T^g]/g!$ is the positive generator of $H^{2g}(J(C),\ZZ)$, we
find:
\begin{eqnarray} 
G_n\cdot\Theta &=&
\sum_{k=0}^{n-1} ({}^{d-g-1}_k)({}^g_{n-k})\frac{(n-k)!}{(n-1-k)!}\\ \nonumber
&=& g\sum_{k=0}^{n-1} ({}^{d-g-1}_k)({}^{g-1}_{n-k-1})\\ \nonumber
&=& g({}^{d-2}_{n-1}).
\end{eqnarray}
Therefore:
\begin{equation} 
[G_k]=({}^{d-2}_{k-1})[C].
\end{equation}

Next we compute the homology class of $\Delta_n$. We use \ref{fies} for
$k=n,\;l=1$: $
\left(\Phi^n_{1,1}\right)_*\Delta_n=\left(\Phi^{n-1}_{1,2}\right)_*H_{n-1}
$.
By induction, similar to \ref{eq5} and using \ref{Gs} (with $n-1$ replaced by
$n$) we have:
\begin{eqnarray}
\left(\Phi^n_{1,1}\right)_*\Delta_n & \ale &
d\left(G_{n-1}-G_{n-2}+\ldots(-1)^{n-3}G_2+(-1)^{n-2}G_1\right)+\\ \nonumber
&&+({}^d_{n-1})2_*C-({}^d_{n-2})3_*C+\ldots+(-1)^{n-1}({}^d_0)(n+1)_*C\\
\nonumber
&\ale &
 d\sum_{k=1}^{n-1}(-1)^{k-1}\mbox{$\frac{1}{k}$}({}^{d-1}_{n-1-k})k_*C+
\sum_{k=2}^{n+1} (-1)^{k} ({}^{d}_{n+1-k})k_*C.
\end{eqnarray}
Taking the homology classes in $J(C)$ we get:
$$
[\Phi^n_{1,1}(\Delta_n)]=\left(d\sum_{k=1}^{n-1}(-1)^{k-1}k({}^{d-1}_{n-1-k})+
\sum_{k=2}^{n+1}(-1)^{k}k^2({}^{d-1}_{n+1-k})\right)[C],
$$
Using \ref{bin3}:
$$
\sum_{k=2}^{n+1}(-1)^k({}^d_{n+1-k})k^2=({}^d_n)+
\sum_{k=0}^{n+1}(-1)^k({}^d_{n+1-k})k^2
=({}^d_n)+({}^{d-3}_{n-1})-({}^{d-3}_{n}),
$$
we find:
$$
[\Phi^n_{1,1}(\Delta_n)]=\left(d({}^{d-3}_{n-2})+({}^d_n)+
({}^{d-3}_{n-1})-({}^{d-3}_n)\right)[C]
$$
and because of
$$
({}^{d-3}_{n-1})-({}^{d-3}_{n})=d({}^{d-3}_{n-1})-(n+1)({}^{d-2}_n), \qquad
({}^{d-2}_{n-1})=({}^{d-3}_{n-1})+({}^{d-3}_{n-2})
$$
we finally obtain:
\begin{equation}\label{homD}
[\Phi^n_{1,1}(\Delta_n)]=\left(d({}^{d-2}_{n-2})-(n+1)({}^{d-2}_n)+({}^d_n)
\right)[C].
\end{equation}

Applying $\left(\Phi_{1,1}^n\right)_*$ to \ref{eq6} and taking homology classes
we get the
following equation for $a,b,c$:
$$
(n+1)({}^{d-2}_{n})=a+b({}^{d-2}_{n-1})+
c\left(d({}^{d-2}_{n-1})-(n+1)({}^{d-2}_n)+({}^d_n)\right).
$$
On the other hand, we have the following intersections numbers in $S_n$:
$$
C\cdot C=G_n\cdot G_n=0, ~~~C\cdot G_n=1,~~~C\cdot \Delta_n=n,~~~ G_n\cdot
\Delta_n =({}^{d-1}_{n-1})
$$
and
$$
H_n\cdot C=d-n,~~~H_n\cdot G_n=({}^{d-1}_n).
$$
Therefore, by intersecting \ref{eq6} with $C$ and $G_n$ respectively, we find
two more equations for $a,b,c$:
$$
d-n=b+cn,~~~~~~~({}^{d-1}_n)=a+c({}^{d-1}_{n-1}),
$$
and \ref{eq2} now follows.
\qed

\subsection{}
With these results it is easy to find many relations between the cycles
$n_*C$.
It is a little surprising that the ones we find are equivalent to $\pi_iC\ale
0$ for $i<2g-d$ and $\pi_{2g-1}C\ale 0$.

\begin{Proposition} \label{recurs}
Let $C$ be a $d$:1 covering of $\PP^1$. Then:
$$
\pi_{2g-i} C \not\ale 0 \quad\Longrightarrow\quad  2\leq i\leq d .
$$
The cycles $\pi_{2g-2}C,\ldots , \pi_{2g-d}C$ span $Z_C/\ale$ and thus
$$
\dim_\QQ\; (Z_C/ \ale) \leq d-1.
$$

For any set $\{n_1,\ldots ,n_{d-1}\}$ of non-zero distinct integers the cycles
$n_{1*}C,\ldots ,n_{d-1*}C$ also span $Z_C/\ale$.

\end{Proposition}
{\bf Proof.}$\quad$
First of all we show that for all $n\in\ZZ$ we have:
\begin{equation}\label{help}
P_nC:=(n+1)_*C-({}^d_{d-1})n_*C+\ldots (-1)^{d}({}^d_0)(n-d+1)_*C+F_d\ale 0
\end{equation}
with $P_n=\G_{n+1}-({}^d_{d-1})\G_n+\ldots\in CH^g(A\times A)$ and with a cycle
$F_d$ depending on $d$ but not on $n$:
$$
F_d=d[-C+({}^{d-1}_{d-3})C^--\mbox{$\frac{1}{2}$}({}^{d-1}_{d-4})2_*C^- +
\ldots (-1)^{d-1}\mbox{$\frac{1}{d-2}$}({}^{d-1}_0)(d-2)_*C^-].
$$
This relation follows from the easily verified identity (for all $n\in\ZZ$):
$$
\left( \Phi^1_{1,n}\right)_* H_1 =\left( \Phi^{d-2}_{-1,n-1}\right)_*H_{d-2}.
$$
Indeed, the l.h.s is by \ref{eq2}:
$$
\left( \Phi^1_{1,n}\right)_* H_1
\ale \left( \Phi^1_{1,n}\right)_* (dC+dG_1-\Delta_1)
\ale dn_*C+dC-(n+1)_*C,
$$
while for the r.h.s. we use
$(\Phi^{k}_{-1,l})_*\Delta_k=(\Phi^{k-1}_{-1,l-1})_*H_{k-1})$:
\begin{eqnarray}
&    & \left( \Phi^{d-2}_{-1,n-1}\right)_*H_{d-2}\\
\nonumber&\ale&
({}^d_{d-2})(n-1)_*C+dG^-_{d-2}-\left(\Phi^{d-3}_{-1,n-2}\right)_*H_{d-3}\\
\nonumber &\ale &
({}^d_{d-2})(n-1)_*C-({}^d_{d-3})(n-2)_*C+\ldots+
(-1)^{d-2}(n-d+1)_*C+ \\ \nonumber
\nonumber& & +d(G_{d-2}^- -G_{d-3}^- +\ldots+(-1)^{d-3}G_1^-)
\end{eqnarray}
(cf. also the proof of \ref{eq3}). From this \ref{help} follows by using
\ref{Gs}.

A more convenient set of relations is obtained from \ref{help} as follows:
\begin{eqnarray}\label{fin}
& & (P_n-P_{n-1})C \; \\ \nonumber
& = & (n+1)_*C-(d+1)n_*C+({}^{d+1}_{d-1})(n-1)_*C+\ldots
+(-1)^{d+1}(n-d)_*C  \\ \label{bz}
& = & \left(\G_{m+d+1}-({}^{d+1}_1)\G_{m+d}+({}^{d+1}_2)\G_{m+d-1}\ldots+
(-1)^{d+1} \G_m\right)C \\ \nonumber
 &  \ale & 0,
\end{eqnarray}
where we substituted: $n:=m+d$.
Relation \ref{bz} can be rewritten as:
\begin{equation}
\G_m*(\G_1-\G_0)^{*(d+1)} C \ale 0\qquad\qquad \forall m\in\ZZ.
\end{equation}

Next we look at the expansion of the $\pi_i$'s in $\G_1-\G_0$:
$$
\pi_i=\mbox{$\frac{1}{(2g-i)!}$}(\log \G_1)^{*(2g-i)}=
$$

$$
=   (\G_1-\G_0)^{*(2g-i)}*
  \left(a_{2g-i}(\G_1-\G_0)^{*0} +\ldots + a_{2g}(\G_1-\G_0)^{*i}\right).
$$
Therefore $\pi_iC\ale 0$ if $2g-i\geq d+1$ so $\pi_{2g-i}C\ale 0$ if $i\geq
d+1$. Since $\pi_{2g}=\Gamma_0$ we also have $\pi_{2g}C=0$, and from \cite{B2},
prop.3 we know $C_{(1)}= 0$ thus:
$$
 C\ale C_{(2)}+C_{(3)}+ \ldots + C_{(d)} \qquad {\rm with}\qquad
 C_{(i)}:=\pi_{2g-i}C
$$
and we conclude that $\dim_\QQ\,( Z_C/\ale) \leq d-1$. We can in fact obtain
$C_{(1)}\ale 0$ from \ref{help}, because a somewhat tedious computation shows
that, for some $c\in \QQ$:
$$
\G_{-1}(\mbox{$\frac{1}{d-1}$}P_0 -\mbox{$\frac{1}{d}$}P_{-1})=c\G_0+\log \G_1
\qquad {\rm mod~} (\G_1-\G_0)^{*(d+1)}* CH^g(A\times A),
$$
here $\G_{-1}\G_{n}=\G_{-n}$ is the product as correspondences, the equality is
in (a quotient of) the ring $CH^g(A\times A)$ with $*$-product. Since $\G_0$
acts trivially on one cycles and $P_nC\ale 0$ the statement follows.
In fact one shows that both sides are equal to, for some $c_1\in\QQ$:
$$
\sum_{n=1}^{d} (-1)^{n-1}\mbox{$\frac{1}{n}$}({}^d_n)\G_n+c_1\G_0.
$$
(Comparing this expression with the one for $u_{d*}G_d$ in \ref{eq3} we see
that
$\pi_{2g-1}C\ale u_{d*}G_d\ale 0$ since $G_d\cong \PP^1$ maps to  a point in
$J(C)$, note that some care must be taken as we defined $G_n$ only for $n<d$).

Finally we observe that since $n_*C_{(i)}=n^iC_{(i)}$, the determinant of the
matrix expressing the $n_{i*}C$ in the $C_{(i)}$ is a Vandermonde determinant,
which is nonzero.
\qed

\begin{Proposition}
\begin{enumerate}
\item For a curve $C$ with $\pi_iC\ale 0$ for $i\leq 2g-4$ (for example any
curve with a $g^1_3$) we have, with $C_{(i)}:=\pi_{2g-i}C$:
\begin{eqnarray}\nonumber
C_{(2)}&\ale & \mbox{$\frac{1}{2}$}(C+C^-)\\ \nonumber
C_{(3)}&\ale  &\mbox{$\frac{1}{2}$}( C-C^-)
\end{eqnarray}
Furthermore: $n_*C\ale \mbox{$\frac{n^3+n^2}{2}$}C\,
-\,\mbox{$\frac{n^3-n^2}{2}$}C^-$.

\item For a curve $C$ with $\pi_iC\ale 0$ for $i\leq 2g-5$ (for example any
curve with a $g^1_4$) we have:
\begin{eqnarray}\nonumber
C_{(2)}&\ale & \mbox{$\frac{-1}{12}$}(2_*C-12C-4C^-)\\ \nonumber
C_{(3)}&\ale  & \mbox{$\frac{1}{2}$}( C-C^-)\\ \nonumber
C_{(4)}&\ale  &\mbox{$\frac{1}{12}$}(2_*C-6C+2C^-).
\end{eqnarray}

\end{enumerate}
\end{Proposition}
{\bf Proof.}$\quad$We give the proof of the second statement, the first being
similar but easier. By assumption we have ($C_{(i)}:=\pi_{2g-i}C$):
$$
C\ale C_{(2)}+C_{(3)}+C_{(4)}\qquad{\rm with}\quad n_*C_{(i)}=n^iC_{(i)}.
$$
Therefore:
\begin{eqnarray} \nonumber
C &\ale & C_{(2)}+C_{(3)}+C_{(4)}\\ \nonumber
C^- &\ale & C_{(2)}-C_{(3)}+C_{(4)}\\ \nonumber
2_*C &\ale & 4C_{(2)}+8C_{(3)}+16C_{(4)}.
\end{eqnarray}
The result follows by straightforward linear algebra.
\qed

\subsection{}
Note that if one specializes a curve $C$ with a $g^1_4$ to a trigonal curve,
then actually $C_{(4)}=\mbox{$\frac{1}{12}$}(2_*C-6C+2C^-)\ale 0$. However, we
could not decide whether for a generic $4$:1 cover of $\PP^1$ we have
$C_{(4)}\ale 0$.

\section{Appendix}
\subsection{}
We recall some facts on binomial coefficients.
The binomial coefficients $({}^n_k)$  are defined (also for negative $n\in \ZZ$
!) as (cf. \cite{ACGH}, VIII.3)
$$
({}^n_k):=\frac{n(n-1) \cdots (n-k+1)}{k(k-1)\cdots 1}\quad (k> 0),\quad
({}^n_0):=1,\quad ({}^n_k):=0\quad (k<0).
$$
With this definition, they are the coefficients in the expansion of $(1+x)^n$:
$$
(1+x)^n=1+({}^n_1)x+\ldots +({}^n_k)x^k+\ldots=\sum_{k=0}^\infty ({}^n_k)x^k.
$$

Comparing the coefficients of $x^k$ in $(1+x)^n(1+x)^m=(1+x)^{n+m}$ one finds,
for all $n,m\in\ZZ$:
\begin{equation}\label{addbin}
\sum_i ({}^n_i)({}^m_{k-i})=(^{n+m}_{k}),
\end{equation}

\begin{Lemma} We have:
\begin{eqnarray}\label{bin1}
\sum_{i=0}^l (-1)^{i}({}^m_{l-i})&=&({}^{m-1}_l),\\ \label{bin2}
\sum_{i=0}^l (-1)^{i-1} i({}^m_{l-i})&=&({}^{m-2}_{l-1}),\\ \label{bin3}
\sum_{i=0}^l (-1)^i i^2({}^m_{l-i})&=&({}^{m-3}_{l-2})-(^{m-3}_{l-1}).
\end{eqnarray}
\end{Lemma}
\vskip0.5truecm
{\bf Proof.}$\quad$Using $({}^{-1}_k)=(-1)^k$ and \ref{addbin}, the first line
can be written as:
$$
\sum_{i=0}^l (-1)^{i}({}^m_{l-i})=\sum_{i=0}^l ({}^{-1}_i)({}^m_{l-i})=
({}^{m-1}_l).
$$
The second line follows in the same way, using $({}^{-2}_{k-1})=(-1)^{k-1}k$:
$$
 \sum_{i=0}^l (-1)^{i-1} i({}^m_{l-i})=
 \sum_{i=0}^l({}^{-2}_{i-1})({}^m_{l-i})=({}^{m-2}_{l-1}).
$$
For the last line we use $({}^{-3}_k)=(-1)^k(k+1)(k+2)/2$, so:
$$
(-1)^kk^2=2({}^{-3}_k)-3({}^{-2}_k)+({}^{-1}_k).
$$
Therefore:
\begin{eqnarray}\nonumber
\sum_{i=0}^l (-1)^i i^2({}^m_{l-i}) &=& 2\sum_{i=0}^l ({}^{-3}_i)({}^m_{l-i})
-3\sum_{i=0}^l ({}^{-2}_i) ({}^m_{l-i})+
\sum_{i=0}^l ({}^{-1}_i) ({}^m_{l-i})\\ \nonumber
&=& 2({}^{m-3}_l)-3({}^{m-2}_l)+({}^{m-1}_l).
\end{eqnarray}
Now use that:
$$
({}^{m-1}_l)=({}^{m-2}_{l-1})+({}^{m-2}_l)=
({}^{m-3}_{l-2})+2({}^{m-3}_{l-1})+({}^{m-3}_l),\quad{\rm and}
$$
$$
({}^{m-2}_l)=({}^{m-3}_{l-1})+({}^{m-3}_l).
$$
\qed

\vspace{\baselineskip}

\noindent Elisabetta Colombo \par\noindent\par\noindent
permanent adress:\par\noindent
Universit\`a di Milano\par\noindent
Dipartimento di Matematica \par\noindent
Via Saldini 51 \par\noindent
20133 Milano\par\noindent
Italia \par\noindent \ \par\noindent
until june 1992:\par\noindent
Department of Mathematics UvA\par\noindent
Plantage Muidergracht 24 \par\noindent
1018 TV Amsterdam\par\noindent
The Netherlands\par\noindent

\vspace{2\baselineskip}

\noindent Bert van Geemen \par\noindent
Department of Mathematics RUU \par\noindent
P.O.Box 80.010 \par\noindent
3508TA Utrecht \par\noindent
The Netherlands \par\noindent

\end{document}